\documentclass[conference]{IEEEtran}
\IEEEoverridecommandlockouts
\usepackage{cite}
\usepackage{booktabs}

\usepackage{pgfplots}
\pgfplotsset{compat=1.16}
\usetikzlibrary{pgfplots.colorbrewer}

\usepackage{amsmath,amssymb,amsfonts,amsthm}
\usepackage{algorithmic}
\usepackage{graphicx}
\usepackage{textcomp}
\usepackage{xcolor}
\def\BibTeX{{\rm B\kern-.05em{\sc i\kern-.025em b}\kern-.08em
    T\kern-.1667em\lower.7ex\hbox{E}\kern-.125emX}}

\newtheorem{theorem}{Theorem}
\newtheorem{lemma}{Lemma}

\newtheorem{corollary}{Corollary}

\usepackage{framed}

\definecolor{applegreen}{rgb}{0.55, 0.71, 0.0}

\usepackage{hyperref}
\hypersetup{
  colorlinks=true,
  citecolor=applegreen
}
\usepackage[nameinlink]{cleveref}

\begin{document}

\title{A Modern Analysis of Hutchinson's Trace Estimator.}

\author{\IEEEauthorblockN{1\textsuperscript{st} Maciej Skorski}
\IEEEauthorblockA{\textit{University of Luxembourg}}}

\maketitle

\begin{abstract}
The paper establishes the new state-of-art in the accuracy analysis of Hutchinson's trace estimator.
Leveraging tools that have not been previously used in this context, particularly 
hypercontractive inequalities and concentration properties of sub-gamma distributions, we offer an elegant and modular analysis, as well as numerically superior bounds. Besides these improvements, this work aims to better popularize the aforementioned techniques within the CS community.
\end{abstract}


\section{Introduction}

\subsection{Background}

\begin{figure}[h!]
\centering
\includegraphics[width=0.75\linewidth]{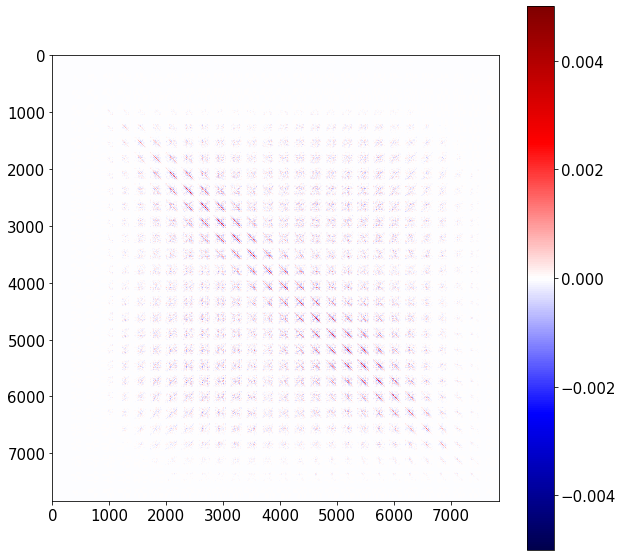}
\caption{Hessian matrix of the linear classifier trained on MNIST data. The diagonal stores valuable information as opposed to rows and columns: entries in rows and colums are nearly zero, whereas the  diagonal average is $\approx 0.25$. The Tensorflow code is available as a Python Notebook at GitHub\protect\footnotemark{}.
}
\label{fig:motivation}
\end{figure}

Estimating the  matrix trace is a problem of fundamental interest~\cite{hutchinson1989stochastic,bai1996some,avron2011randomized}
and arises in many problems such as approximating spectral properties of matrices~\cite{lin2016approximating,han2017approximating,di2016efficient},
solving partial differential equations~\cite{haber2012effective,van2012adaptive,young2012application,roosta2014stochastic,van2011seismic}, error evaluation in machine-learning~\cite{golub1997generalized}, and combinatorial counting~\cite{avron2010counting} to mention a few. For readers particularly interested in data science or optimization, it is of critical interest as the hessian trace stores valuable information about curvature. To give a meaningful example, consider fitting a linear classifier on the MNIST dataset~\cite{lecun1998gradient} of hand-written digits, widely used as a benchmark and a toy problem in machine learning. Since images are in resolution $28\times 28$ and grouped in $10$ classes, the number of parameters is $m=28\cdot28\cdot 10 = 7840$, and the size of the hessian matrix is $m^2 \approx 6\cdot 10^6$.
The diagonal average, proportional to the trace, equals the average eigenvalue (valuable information), whereas individual rows or columns are nearly zero (up to random noise). This is illustrated in \Cref{fig:motivation}.

\footnotetext{\url{https://github.com/maciejskorski/ml_examples/blob/master/TF_Hess_Visualize.ipynb}}

As seen from the above example, already toy examples lead to huge matrices; for larger problems storing and inspecting such matrices is impossible. Fortunately, it is possible to estimate the trace (equal to the diagonal sum, or equivalently to the sum of eingenvalues) \emph{without knowing the full matrix}. The  requirement is that one can compute efficiently \emph{matrix-vector} products; in the case of the hessian such products can be computed by by auto-differentiation~\cite{christianson1992automatic} in all popular machine-learning frameworks such as Tensorflow~\cite{abadi2016tensorflow}, PyTorch~\cite{paszke2019pytorch} and JAX~\cite{jax2018github}. Under the assumption that an $m\times m$ matrix $A$ is 
\emph{symmetric positive semi-definite} (which applies to all hessians), the popular estimator due to Hutchinson~\cite{hutchinson1989stochastic} is 
\begin{align}\label{eq:h_estimator}
 \mathrm{tr}_H(A)\triangleq  z^T A z,\quad z\sim^{iid} \{-1,1\}^{m},
\end{align}
where the components of the vector $z$ are \emph{Rademacher random variables}, equal to $+1$ or $-1$ with equal probability. The above estimator uses only one sampled $v$ and thus is noisy (of high-variance), so it is usually boosted by averaging $n$ independent trials (for suitably large $n$). Formally:
\begin{align}\label{eq:h_estimator_n}
 \mathrm{tr}_{H^{(n)}}(A)\triangleq  \frac{1}{n}\sum_{i=1}^{n}z_i^T A z_i,\quad z_i\sim^{iid} \{-1,1\}^{m}.
\end{align}
The estimators are unbiased, that is correct on average:
\begin{align}\label{eq:h_estimator_unbiased}
\mathbb{E} \mathrm{tr}_H(A) = \mathbb{E} \mathrm{tr}_{H^{(n)}}(A) = \mathrm{tr}(A),
\end{align}
and this is fairly easy to prove. The focus of this work is on a more challenging problem of
understanding their \emph{concentration properties}. Here we ask the following question:
\begin{framed}
\begin{quote}
    What sample size $n$ \emph{guarantees} the desired relative error $\epsilon$ at the confidence level of $1-\delta$? 
\end{quote}
\end{framed}
Formally, for the error $\mathrm{err}_{H^{(n)}}(A)\triangleq \frac{\mathrm{tr}_{H^{(n)}}(A)}{  \mathrm{tr}(A)}-1$, and fixed $\epsilon,\delta$, we want possibly small $n=n(\epsilon,\delta)$ such that
\begin{align}\
 | \mathrm{err}_{H^{(n)}}(A)|  \leqslant \epsilon \quad \text{with prob. at least}\ 1-\delta.
\end{align}

\subsection{Related Work}

The estimator is already more than 30 years old~\cite{hutchinson1989stochastic}, and although  alternatives exist (such as methods based on Lanczos quadrature~\cite{ubaru2017fast}), it is provably best in terms of variance~\cite{roosta2015improved} and arguably wins in simplicity, being preferred by developers of statistical/learning software~\cite{yao2019pyhessian}. Quite surprisingly, it had been used for a while without a rigorous assessment of accuracy, until the work of Avron and Toledo~\cite{avron2011randomized}, who established the finite sample size guarantees. Their result was later improved by Roosta-Khorasani and Ascher~\cite{roosta2015improved}, who essentially got rid of the lossy proof step involving a crude union bound. Their approach is based on the Chernoff-like direct calculations,  and offers the bound $n=O(\epsilon^{-2}
\log(1/\delta))$.

\section{Contribution}

\subsection{Summary}
In this work we offer a \emph{modular, modern and more accurate} analysis of Hutchinson's estimator. The improvements upon prior works can be summarized as follows:
\begin{itemize}
    \item \textbf{Modularity and Novelty of Techniques.} In the first step, we explicitly link the estimator accuracy to the dispersion of the \emph{Rademacher Quadratic Chaos} with a unitary matrix kernel. We then obtain a bound on such chaoses by means of \emph{Hypercontractive Inequality},
    an important result in Boolean Fourier Analysis~\cite{o2014analysis}, originally due to Aline Bonami~\cite{Bonami1968}. As the final step we express this bound as the \emph{sub-gamma property}  (see the works of Boucheron~\cite{boucheron2013concentration}), which makes it convenient (and accurate) to conclude the final concentration result.
    Moreover, we state our results for the multiple-sample \emph{and} one-shot estimators; this is of interest, since the base building block may be boosted differently than by averaging (see for example the median algorithm~\cite{vershynin2018high}).
    
    Thus, in contrast to ad hoc calculations in prior works, we are able to use established tools from the field of high-dimensional probability and Fourier analysis. Our transparent approach opens a way for further refinements.
    \item \textbf{Superior Accuracy.} With the dedicated tools we obtain bounds that are numerically better up to an order of magnitude. They are stated as elementary formulas, convenient to use in  practical applications of trace estimation.
\end{itemize}

\subsection{Preliminaries}

Below we explain the notation and definitions used when formulating our results presented in the next section.

The \emph{$d$-th norm} of a r.v. $X$ is $\|X\|_d \triangleq (\mathbb{E}|X|^d)^{1/d}$; it is a valid norm (e.g. it obeys the triangle inequality) when $d\geqslant 1$~\cite{vershynin2018high}. A useful tool for studying  concentration is the \emph{moment generating function}, defined as $
\mathrm{MGF}_X(t) = \mathbb{E}\exp(tX)$ (a function of real parameter $t$). In modern high-dimensional probability one classifies distributions according to the behaviour of $\mathrm{MGF}$; in particular we call $X$ \emph{sub-gamma} with variance factor $v$ and scale $c$ when $\log\mathrm{MGF}_{|X|}(t) \leqslant \frac{v t^2}{2(1-c t)}$, and denote by $X\in\Gamma(v,c)$~\cite{boucheron2013concentration};
the same formula holds for the centered gamma distribution, hence the name.

\subsection{Results}

In what follows we assume that $A$ is non-zero symmetric positive semi-definite matrix. Then, necessarily, $\mathrm{tr}(A)>0$.

\subsubsection{One-Shot Estimator}

Below we present the following bound on the accuracy of the estimator in~\eqref{eq:h_estimator}:

\begin{theorem}[Hutchison's Estimator 1-Sample Bound]\label{thm:main}
For any integer $d\geqslant 2$, the relative error of Hutchison's Estimator~\eqref{eq:h_estimator} obeys the following bound
\begin{align}\label{eq:one_shot_moment}
\| \mathrm{err}_H(A) \|_d \leqslant d-1.
\end{align}
In particular, this implies the sub-gamma behaviour
\begin{align}\label{eq:sub_gamma_1}
\mathrm{err}_H(A) \in \Gamma\left(1,\frac{8}{3}\right),
\end{align}
which in turns gives the following probability tail bound
\begin{align}\label{eq:tail_1}
\Pr[ |\mathrm{err}_{H}(A)|\geqslant \epsilon]\leqslant \mathrm{e}^{-\frac{\epsilon^2}{2\left(1-\frac{8}{3}\epsilon\right)}}.
\end{align}
\end{theorem}
The proof goes along the
following lines: a) we express the error in form of Rademacher chaos b) we use hypercontractive inequalities to bound its moment; this establishes the part \eqref{eq:one_shot_moment}
c) we plug the moment bound into Taylor's expansion of the moment generating function and then estimate the series, arriving at the sub-gamma condition \eqref{eq:sub_gamma_1}.
By the tail properties of sub-gamma distributions, we conclude the tail bound \eqref{eq:tail_1}.

\subsubsection{Multiple-Sample Estimator}

Next, we move to the multiple-sample estimator in \eqref{eq:h_estimator_n}. We establish the following

\begin{theorem}[Hutchinson's Estimator $n$-Sample Bound]\label{thm:main2}
The relative error of $n$-sample Hutchinson's Estimator has the following sub-gamma behaviour
\begin{align}
\mathrm{err}_H(A) \in \Gamma\left(\frac{1}{n},\frac{8}{3}\right).
\end{align}
In particular, this implies the following bound on the error tail probability, valid for any $ 0<\epsilon<3/8$:
\begin{align}
\Pr[ |\mathrm{err}^{(n)}_{H}(A)|\geqslant \epsilon]\leqslant \mathrm{e}^{-\frac{n\epsilon^2}{2\left(1-\frac{8}{3}\epsilon\right)}}.
\end{align}
\end{theorem}
This result is obtained from \Cref{thm:main}, by using extra facts on the concentration of sums of sub-gamma random variables: essentially the variance factor decreases form $1$ to $\frac{1}{n}$ for sums of length $n$ (just as the variance of $n$ iid terms).

\begin{corollary}[Sample Size for Hutchinson's Estimator]\label{cor:sample_size}
The relative error is absolutely bounded by $\epsilon$ with probability $1-\delta$, provided that the sample size $n$ is bigger or equal to
\begin{align}
    n(\epsilon,\delta) = \frac{2\left(1-\frac{8}{3}\epsilon\right)\log\left(\frac{1}{\delta}\right)}{\epsilon^2}.
\end{align}
This holds for any $0<\delta<1$ and $0<\epsilon<3/8$.
\end{corollary}

\newpage

\subsection{Numerical Comparison with Related Work}

Below we demonstrate numerical improvements with respect to the prior works
~\cite{avron2011randomized,roosta2015improved}. The bounds, with exact constants, are summarized in \Cref{tab:my_label} below. 

\begin{table}[h!]
    \centering
    \begin{tabular}{c|c}
    \toprule
      Author   &  Sample Size $n=n(\epsilon,\delta)$ \\  
      \midrule
      \textbf{this work}  & $ \frac{2\left(1-\frac{8}{3}\epsilon\right)\log\left(\frac{1}{\delta}\right)}{\epsilon^2}$ \\[1.25ex]
       
      \cite{roosta2015improved} & $\frac{6\log\left(\frac{2}{\delta}\right)}{\epsilon^2}$  \\[1.25ex] 
      \cite{avron2011randomized} & $\frac{6\log\left(\frac{2}{\delta}\right)+6\mathrm{rank}(A)}{\epsilon^2}$\\[1.25ex] \bottomrule
    \end{tabular}
     \bigskip
    \caption{Bounds on the accuracy of Hutchinson's Estimator.}
    \label{tab:my_label}
\end{table}

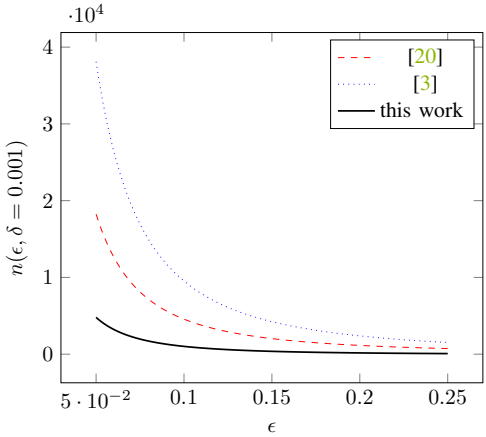
\begin{figure}[h!]
\centering
\resizebox{.75\linewidth}{!}{
\begin{tikzpicture}
[
  declare function = {
    sample_my(\eps,\de) = (2*(1-8*\eps/3)*-ln(\de))/(\eps^2);
    sample_1(\eps,\de,\rnk) = (6*ln(1/\de)+6*ln(\rnk))/(\eps^2);
    sample_2(\eps,\de) = (6*ln(2/\de))/(\eps^2);
  }
]
 \begin{axis}[xlabel={$\epsilon$},ylabel={$n(\epsilon,\delta=0.001)$},cycle list name=color list],
  \addplot+[domain=0.05:0.25, samples=200, variable=\eps, dashed] plot ({\eps}, {sample_2(\eps,0.001)});
 \addplot+[domain=0.05:0.25, samples=200, variable=\eps, dotted] plot ({\eps}, {sample_1(\eps,0.001,7840)});
  \addplot+[domain=0.05:0.25, samples=200, variable=\eps, thick] plot ({\eps}, {sample_my(\eps,0.001)});
  \legend{\cite{roosta2015improved}, \cite{avron2011randomized}, {this work}};
 \end{axis}
\end{tikzpicture}
}
\caption{Sample size for Hutchinson's Estimator. We assume $
\delta=10^{-3}$ which results in confidence of $0.999$. The matrix size is $7840$, which corresponds to the Hessian of the linear MNIST classification model (input data are black-white images in the resolution of $28\times 28$, grouped in $10$ classes).}
\label{fig:sample_size}
\end{figure}

Furthermore, in~\Cref{fig:sample_size}, we present a detailed evaluation of how the sample size $n$ depends on the error $\epsilon$, when the confidence parameter $\delta$ is fixed. The setup for this experiment is the discussed hessian of a toy classifier on MNIST.

\section{Auxiliary Results}

\subsubsection{Matrix Theory}

We will need the following fact on the decomposition of symmetric matrices~\cite{mello2017invitation}.
Recall that a matrix $B$ is orthonormal when $B^T = B^{-1}$ ($\cdot^T$ denotes transposition); in particular,  the columns and rows of $B$ are of unit length.
\begin{lemma}[Symmetric Matrix Factorization]\label{lemma:matrix_factorize}
Any symmetric real matrix $A$ can be written as $A = B^T \Lambda B$ where $B$ is an orthonormal real matrix and $\Lambda$ is diagonal, consisting of (necessarily real) eigenvalues of $A$.
\end{lemma}

\subsubsection{Quadratic Chaos}

We will need the following special case of the celebrated 
Hypercontractivity Inequality~\cite{Bonami1968,o2014analysis}:

\begin{lemma}[($2,d)$-Hypercontractivity]\label{lemma:Hypercontractivity}
Let $F$ be a polynomial of degree $2$ in Rademacher variables $Z_i$ (e.g. $F=\sum_{i\not=j}a_{i,j}Z_{i}Z_{j}$ for some weights $a_{i,j}$). Then
\begin{align}
   \|F\|_d  \leqslant (d-1)\|F\|_2.
\end{align}
\end{lemma}

\subsubsection{Sub-Gamma Random Variables}

The tail probability of sub-gamma distributions can be estimated from the MGF bound by the Crammer-Chernoff method, as follows:
\begin{lemma}[Sub-Gamma Tails]\label{lemma:subgamma_tail}
Let $ X\in \Gamma(v,c)$, then
\begin{align}
 \forall \epsilon>0:\ \Pr[|X|\geqslant \epsilon] \leqslant \mathrm{e}^{-\frac{\epsilon^2}{v+c\epsilon}}.
\end{align}
\end{lemma}
Moreover, for sub-gamma distributions the following composition property holds for sums of independent terms:
\begin{lemma}[Sub-Gamma Sums]\label{lemma:subgamma_sum}
Let random variables $ X_i\in \Gamma(v_i,c_i)$ be independent, then
\begin{align}
    \sum_i X_i \in \Gamma(v',c'),\ v'\triangleq \sum_i v_i,c'\triangleq \max_i c_i.
\end{align}
\end{lemma}
The proofs of these lemmas can be found in \cite{boucheron2013concentration}.

\section{Proofs}

\subsection{Proof of \Cref{thm:main}}

Define the following random variable
\begin{align}\label{eq:offset}
Y = z^T A z - \mathbb{E}[ z^T A z ] = z^T A z - \mathrm{tr}(A).
\end{align}
then our task is to bound the concentration of $Y$.

\subsubsection{Reduction to Off-Diagonal Quadratic Chaos}

Since $A$ is symmetric, by \Cref{lemma:matrix_factorize} we have $A = B^T \Lambda B$ where $B$ is orthonormal and $\Lambda$ is diagonal made of the eigenvalues of $A$. Since, in addition, the matrix $A$ is positive semi-definite, its eigenvalues are non-negative. Thus, we can write
\begin{align}
z^T A z = (Bz)^T \Lambda (Bz) = \|\Lambda^{1/2} Bz\|_2^2.
\end{align}
This can be decomposed as
\begin{align}
z^T A z  = \sum_i Y_i    ,\quad Y_i & \triangleq (\Lambda^{1/2} Bz)_i^2.
\end{align}
To further simplify, denote by $\lambda_1,\ldots,\lambda_m$ all the eigenvalues of $A$; these are precisely the diagonal entries of $\Lambda$. We note that $Y_i = \lambda_i (\sum_{j}B_{i,j}z_j)^2$,
and $\mathbb{E}Y_i =\lambda_i \sum_{j} B_{i,j}^2 $ (because $\mathbb{E}z_j^2=1$). Therefore we obtain
\begin{align}\label{eq:quad_chaos}
    Y_i-\mathbb{E}Y_i =\lambda_i \sum_{j\not=j'} B_{i,j} B_{i,j'} z_j z_{j'}.
\end{align}
Thus, every $Y_i$ is a quadratic form in $z_j$ with zero-diagonal.

\subsubsection{Bounding Quadratic Chaos}

For every fixed $i$ we apply \Cref{lemma:Hypercontractivity} to
$Y_i-\mathbb{E}Y_i$, which is a polynomial of degree $2$ in Rademacher random variables $z_i$ (note 
that the off-diagonal property guarantees the degree is exactly two).

Since the random variables $z_{j}z_{j'}$ indexed by tuples $(j,j')$ for $j\not=j'$ are uncorrelated, 
we have $\|Y_i-\mathbb{E}Y_i\|_2^2=\lambda_i^2\sum_{j\not=j'}B_{i,j}^2B_{i,j'}^2$;
but $B$ is orthonormal, so $\sum_{i,j}B_{i,j}^2=1$ and $\|Y_i-\mathbb{E}Y_i\|_2^2\leqslant \lambda_i^2$. Thus, we obtain
\begin{align}
    \|Y_i-\mathbb{E}Y_i\|_d \leqslant \lambda_i^2(d-1).
\end{align}
By the triangle inequality 
\begin{align}
    \|\sum_i (Y-\mathbb{E}Y_i)\|_d \leqslant (d-1)\sum_i \lambda_i.
\end{align}
Finally by the definition of $Y_i$ and the standard identity $\mathrm{tr}(A)=\sum_i\lambda_i$ from linear algebra~\cite{mello2017invitation}
\begin{align}\label{eq:moment_bound}
    \| z^T A z - \mathrm{tr}(A) \|_d \leqslant (d-1)\mathrm{tr}(A).
\end{align}
Dividing the both sides of this inequality by $\mathrm{tr}(A)$, we complete the proof of the first part of \Cref{thm:main}.

\subsubsection{Bounding Moment Generating Function}

Let $E = z^T A z / \mathrm{tr}(A)-1 $ be the relative error of the estimator. Then by the previous result, the Taylor expansion $\mathrm{e}^{tx}=\sum_{d\geqslant 0} (tx)^d/d!$, and the fact that $\mathbb{E}E=0$ we have
\begin{align}
  MGF(E) \leqslant 1+ \sum_{d\geqslant 2}\frac{ t^d(d-1)^d}{d!}.
\end{align}
Let $a_d=(d-1)^d/d!$. Then we have
\begin{align}
    \frac{a_{d+1}}{a_d} = \frac{d^{d+1}}{(d-1)^{d}} \cdot \frac{1}{d+1}.
\end{align}
This ratio is decreasing when $d\geqslant 2$, and thus
\begin{align}
    \frac{a_{d+1}}{a_d} \leqslant \frac{8}{3}.
\end{align}
Therefore, it follows that
\begin{align}
    \mathrm{MGF}(E) \leqslant 1+\frac{t^2}{2(1-\frac{8}{3}t)},
\end{align}
which, by the standard inequality $\log(1+u)\leqslant u$, implies
\begin{align}
    \log\mathrm{MGF}(E) \leqslant \frac{t^2}{2(1-\frac{8}{3}t)},
\end{align}
so that, by definition, $E\in\Gamma(1,8/3)$. This completes the proof of the second part of \Cref{thm:main}.

\subsubsection{Concluding Concentration Properties}

The third part of \Cref{thm:main} follows directly by \Cref{lemma:subgamma_tail}.

\subsection{Proof of \Cref{thm:main2}}

We first apply the result on sum of independent sub-gamma distributions from \Cref{lemma:subgamma_sum}, which proves the first part of the theorem. The second part follows by the result on sub-gamma tails, stated in  \Cref{lemma:subgamma_tail}.

\subsection{Proof of \Cref{cor:sample_size}}

The corollary follows by equating the right-hand side of the tail bound from \Cref{thm:main2}, and equating with $\delta$. Solving with respect to $n$ gives the claimed formula on the sample size.

\section{Conclusion}

This work establishes a new state-of-art bound on the classical trace estimator due to Hutchinson.
The core idea is the clever usage of bounds on Rademacher Chaos and the theory of sub-gamma distributions. The results are numerically superior and of immediate interest to any problems involving trace estimation. Besides, the author hopes to contribute to better popularization of the discussed techniques, novel in the problem context, within the computer science community.

\bibliographystyle{IEEEtran}
\bibliography{citations}

\end{document}